\documentclass[12pt,twoside,a4paper]{article}
\usepackage{times}
\usepackage{bm}
\usepackage{braket}
\usepackage{graphicx}
\usepackage{theorem}
\usepackage{amsmath}
\usepackage{amssymb,mathrsfs}
\usepackage{latexsym}
\usepackage{natbib}
\usepackage{algorithm}
\usepackage{algpseudocode}
\usepackage{caption}
\newfloat{procedure}{H}{loa}
\floatname{procedure}{Procedure} 
\usepackage[flushmargin]{footmisc}

\theorembodyfont{\itshape}
\newtheorem{theorem}{Theorem}[section]
\newtheorem{lemma}{Lemma}[section]

\theorembodyfont{\rmfamily}
\newtheorem{definition}{Definition}[section]{\bf}{\rm}
\newtheorem{assumption}{Assumption}[section]{\bf}{\rm}

\newtheorem{remark}{Remark}[section]{\itshape}{\rmfamily}
{\itshape}{\rmfamily}
\newenvironment{proof}{\noindent{\it Proof.~~}}{\medskip}

\makeatletter
\def\eqnarray{\stepcounter{equation}\let\@currentlabel=\theequation
\global\@eqnswtrue
\global\@eqcnt\z@\tabskip\@centering\let\\=\@eqncr
$$\halign to \displaywidth\bgroup\@eqnsel\hskip\@centering
  $\displaystyle\tabskip\z@{##}$&\global\@eqcnt\@ne 
  \hfil$\;{##}\;$\hfil
  &\global\@eqcnt\tw@ $\displaystyle\tabskip\z@{##}$\hfil 
   \tabskip\@centering&\llap{##}\tabskip\z@\cr}
\makeatother
\makeatletter
    \renewcommand{\theequation}{%
    \thesection.\arabic{equation}}
    \@addtoreset{equation}{section}
  \makeatother
\makeatletter
\def\Left#1#2\Right{\begingroup%
   \def\ts@r{\nulldelimiterspace=0pt \mathsurround=0pt}%
   \let\@hat=#1%
   \def\sht@im{#2}%
   \def\@t{{\mathchoice{\def\@fen{\displaystyle}\k@fel}%
          {\def\@fen{\textstyle}\k@fel}%
          {\def\@fen{\scriptstyle}\k@fel}%
          {\def\@fen{\scriptscriptstyle}\k@fel}}}%
   \def\g@rin{\ts@r\left\@hat\vphantom{\sht@im}\right.}%
   \def\k@fel{\setbox0=\hbox{$\@fen\g@rin$}\hbox{%
      $\@fen \kern.3875\wd0 \copy0 \kern-.3875\wd0%
      \llap{\copy0}\kern.3875\wd0$}}%
      \def\pt@h{\mathopen\@t}\pt@h\sht@im%
      \Right}%
\def\Right#1{\let\@hat=#1%
   \def\st@m{\mathclose\@t}%
   \st@m\endgroup}
\makeatother

\newcommand{\vc}{\bm}
\makeatletter
\DeclareRobustCommand\widecheck[1]{{\mathpalette\@widecheck{#1}}}
\def\@widecheck#1#2{%
    \setbox\z@\hbox{\m@th$#1#2$}%
    \setbox\tw@\hbox{\m@th$#1%
       \widehat{%
          \vrule\@width\z@\@height\ht\z@
          \vrule\@height\z@\@width\wd\z@}$}%
    \dp\tw@-\ht\z@
    \@tempdima\ht\z@ \advance\@tempdima2\ht\tw@ \divide\@tempdima\thr@@
    \setbox\tw@\hbox{%
       \raise\@tempdima\hbox{\scalebox{1}[-1]{\lower\@tempdima\box
\tw@}}}%
    {\ooalign{\box\tw@ \cr \box\z@}}}
\makeatother

\newcommand{\ol}{\overline}

\newcommand{\vertiii}[1]%
{{\left\vert\kern-0.25ex\left\vert\kern-0.25ex\left\vert #1 
 \right\vert\kern-0.25ex\right\vert\kern-0.25ex\right\vert}}
\newcommand{\down}[2]{\smash{\lower#1\hbox{#2}}}
\newcommand{\up}[2]{\smash{\lower-#1\hbox{#2}}}

\newcommand{\dm}{\displaystyle}

\newcommand{\EE}{\mathsf{E}}
\newcommand{\PP}{\mathsf{P}}

\newcommand{\one}{\mbox{$1$}\hspace{-0.25em}{\rm l}}

\newcommand{\calA}{\mathcal{A}}

\newcommand{\calC}{\mathcal{C}}

\newcommand{\calN}{\mathcal{N}}

\newcommand{\calP}{\mathcal{P}}

\newcommand{\bbA}{\mathbb{A}}

\newcommand{\bbN}{\mathbb{N}}
\newcommand{\bbR}{\mathbb{R}}

\newcommand{\bbT}{\mathbb{T}}

\newcommand{\bbZ}{\mathbb{Z}}
\DeclareMathOperator*{\argmin}{arg\,min}

\newcommand{\dd}[1]{\if#11 1\!\!1 
\else {\if#1C I\!\!\!C
\else {\if#1G I\!\!\!G 
\else {\if#1J J\!\!\!J 
\else {\if#1S S\!\!\!S
\else {\if#1Z Z\!\!\!Z
\else {\if#1Q O\!\!\!\!Q
\else I\!\!#1
\fi}
\fi}
\fi}
\fi} 
\fi} 
\fi} 
\fi} 

\makeatother

\begin{document}\thispagestyle{empty} 

\phantom{}\vspace{-10mm}

\hfill

{\Large{\bf
\begin{center}
A discrete-time single-server Poisson queueing game: Equilibria simulated by an agent-based model
%
%
%
%
\end{center}
}
}

\begin{center}
\textbf{Yutaka Sakuma${}^{a}$%
\footnotetext[1]{
Department of Computer Science, National Defense Academy of Japan, Kanagawa 239-8686, Japan
\\ 
(E-mail: sakuma@nda.ac.jp)
},
Hiroyuki Masuyama${}^{b}$%
\footnotetext[2]{
Department of Systems Science, Graduate School of Informatics, Kyoto University, Kyoto 606-8501, Japan
\\
(E-mail: masuyama@sys.i.kyoto-u.ac.jp)}, and 
Emiko Fukuda${}^{c}$\footnotetext[3]{
Department of Industrial Engineering and Economics, Tokyo Institute of Technology, Tokyo 152-8550, Japan
\\
(E-mail: fukuda.e.ac@m.titech.ac.jp)}
}

\bigskip
\medskip

{\small
\textbf{Abstract}

\medskip

\begin{tabular}{p{0.85\textwidth}}
%
This paper considers a discrete-time single-server queue with a single acceptance period for a Poissonian population of homogeneous customers. Customers are served on a first-come first-served (FCFS) basis, and their service times  are independent and identically distributed with a general distribution. We assume that each customer chooses her/his arrival-time slot with the goal of minimizing her/his expected waiting time in competition with other customers. For this queueing game, we derive a symmetric (mixed-strategy) Nash equilibrium; that is, an {\it equilibrium arrival-time distribution} of homogeneous customers, where their expected waiting times are identical. We also propose an agent-based model, which simulates the dynamics of customers who try to minimize their waiting times for service. Through numerical experiments, we confirm that this agent-based model achieves, in steady state, an arrival-time distribution similar to the equilibrium arrival-time distribution analytically obtained.
\end{tabular}
}
\end{center}

\begin{center}
\begin{tabular}{p{0.90\textwidth}}
{\small

\medskip

}
\end{tabular}

\end{center}

\section{Introduction}
\label{sec:Introduction}
Many real-life service systems have acceptance periods for arriving customers, for  example,  hospitals, call centers, and check-in counters in airports, etc. In such a system, the number of arriving customers is not specified and usually large, and each customer faces with the decision problem of when to arrive at the system so as to achieve a certain goal,  typically, to minimize the waiting time for service in competition with other customers. Such queueing models are called {\it (non-cooperative) queueing games}.

\cite{Glaz83} pioneered the study on queueing games with opening and closing times, and an unspecified number of arriving customers. Their queueing game is defined in continuous time. The system has a single server with an infinite buffer. The system is open for a Poissonian population of homogeneous customers, though it has an acceptance period for arriving customers, which starts from the opening time and ends at the closing time. Note that {\it early arrival} is allowed therein.  Thus, customers can arrive at the system before the opening time, but even though they take this option, serving them is postponed until the opening time. In summary, customers choose their arrival times which are before the closing time. Customers also demand service times independently and identically distributed (i.i.d.)\ according to an exponential distribution, and their demands are processed by the server on a first-come first-served (FCFS) basis. In this setting, each customer chooses her/his arrival time with the goal of minimizing her/his waiting time in competition with other customers. This queueing game is a non-cooperative one with a random number of customers. For the queueing game, \cite{Glaz83} derive an equilibrium strategy of arriving customers; that is, an arrival-time distribution such that the expected waiting times of customers are identical. For convenience, we refer to such an arrival-time distribution as an {\it equilibrium arrival-time distribution}, and denote by an {\it equilibrium expected waiting time}, an expected waiting time of an arbitrary customer with an equilibrium arrival-time distribution.

There are several studies on continuous-time queueing games, following \cite{Glaz83}'s work. We provide a brief survey (for a comprehensive one, see \citealt{Hass16}).  \cite{Hass11} obtain an equilibrium arrival-time distribution in a queueing game without early arrivals. \cite{Hass11} also report that no early arrivals reduce the equilibrium expected waiting time, especially when the system is heavily loaded (see Figure 3 therein). In addition to waiting cost, \cite{Havi13} and \cite{Ravn14} consider queueing games taking into account tardiness cost and arrival order cost, respectively. For a queueing game with tardiness cost, \cite{June13} have proved the existence and uniqueness of the equilibrium arrival-time distribution, under the assumptions of early arrivals, no closing time, exponential service times, and a random but bounded number of customers. 

The studies mentioned above assume that service times are exponentially distributed. Unlike these studies, \cite{Brei17} extends the queueing game of \cite{June13} to a case where more general customers' preferences and service-time distribution are allowed (though early arrivals are not allowed). However, in developing the numerical procedure for an equilibrium arrival-time distribution, \cite{Brei17} assumes that service times are exponentially distributed (see Lemma 8 therein).

In this paper, we consider a queueing game with a general service-time distribution, though it is in discrete time. \cite{Ostu08} propose a discrete version of \citealt{Vick69}'s transportation models (which can be considered queueing games with a fixed number of players demanding deterministic service times), for which they establish an algorithm for computing equilibrium. In the similar vein, we consider a discrete-time queueing game, which makes it easier to handle general service-time distributions.  However, even for discrete-time models (of course, for continuous-time ones), there are no previous studies that consistently assume that the service-time distribution is general. Although \cite{Rapo04} and \cite{Seal05} study discrete-time queueing games, they assume that the number and service times of customers are fixed to be constants. 

As mentioned above, a feature of our model is not to specify the type of the service-time distribution. More specifically, our model is a discrete-time FCFS single-server queueing game with an acceptance period for the Poisson population of customers having a general service-time distribution. For simplicity, we refer to such a queueing game as a {\it (discrete-time) single-server Poisson queueing game}. We also assume that no early arrivals are allowed in this paper.

For this single-server Poisson queueing game, we establish an algorithm for computing an equilibrium arrival-time distribution. Through some numerical examples, we also show that the large variation of service times causes the rush of customers to the opening slot. In addition, we propose an agent-based model simulating the behavior of non-cooperative customers who try to minimize their waiting times according to their experiences, i.e., their own histories of waiting times that are accumulated until their respective decision times. By running this agent-based model for a sufficiently long time, we have an arrival-time distribution similar to the equilibrium arrival-time distribution computed by our algorithm. This result implies that our queueing game can serve as a mathematical model of the real behavior of non-cooperative customers.

The rest of the paper is organized as follows. Section~\ref{sec:Model} describes our queueing game and its equilibrium arrival-time distribution. Section~\ref{sec:Expected waiting time} provides fundamental results on the expected unfinished workload and waiting time. Section~\ref{sec-arrival-distribution} establishes an algorithm for computing an equilibrium arrival-time distribution, and presents some numerical examples. Section~\ref{sec:Simulation_Model} proposes an agent-based model. Section~\ref{sec:Conclusion} is devoted to concluding remarks.

Finally, we provide the notation and basic definitions used in the subsequent sections.
Let $\bbN = \{1, 2, 3, \ldots\}$ and $\bbZ_{+} = \{0\} \cup \bbN$, and let $\bbR_{+}$ denote the set of all nonnegative real numbers. 
Let $\one(\cdot)$ denote the indicator function of the event in the parenthesis.
For any set $\calA$, let $|\calA|$ denote its cardinality. Furthermore, the empty sum is defined as zero, e.g., $\sum_{n=1}^{0} ( \cdot ) = 0$.
\section{Model description}
\label{sec:Model}

We consider a discrete-time single-server queue with an infinite waiting room. The time axis is divided into (time) slots of length one, and thus the $k$-th slot ($k\in\bbZ_+$) is the interval $[k,k+1]$, which is referred to as slot $k$. The system is empty immediately before slot 0 starts, and it accepts customers arriving in slots 0 through $T$ ($T \in \bbN$), which are the acceptance period. Accepted customers are served on the first-come first-served (FCFS) basis, though the customers arriving in the same slot are served in random order. Note here that, even after the acceptance period ends, the server continues to work until the queue is empty.

For the discrete-time queueing model, we have to make a condition on which instants in slots the events (such as arrivals, departures, and starts of services) occur.
\begin{assumption}[Instants of arrivals, departures, and starts of services]
\label{assumpt-V_0-}
In each slot $t \in \bbT:=\{0,1,\dots,T\}$, the arrival of customers can occur immediately after the slot starts, i.e., at time $t+ := \lim_{\epsilon \downarrow 0} (t + \epsilon)$. In addition, if the unfinished workload (in the system), which is the total of the remaining service demand of customers in the system, is positive in the middle of slot $t$ ($t \in \bbZ_+$), then the server handles and decreases it by one immediately before the end of the slot; that is, the unfinished workload (if any) decreases by one at time $(t+1)- := \lim_{\epsilon \downarrow 0}(t+1-\epsilon)$. Thus, the departure of customers can occur immediately before each slot ends. 
\end{assumption}
\begin{assumption}[Service time distribution]
\label{assumpt-service}
The service times of customers are i.i.d.\ with discrete distribution $(b(k);k\in\bbN)$ having finite positive mean $\beta:=\sum_{k=1}^{\infty} k b(k)$, where $b(0) = 0$ is assumed.
\end{assumption}

\begin{remark} 
Consider a continuous-time queueing model with acceptance period $[0,T]$ and service-time distribution $G:\bbR_+\to[0,1)$, where $G(0) = 0$. Divide the time axis into slots of length $d:=T/M$, where $M>0$ is a sufficiently large integer. Let $b(k) = G(kd) - G((k-1)d)$ for $k \in \bbN$. This discretization yields a discrete-time model as an approximation to the original continuous-time queueing model.
\end{remark}

\begin{figure}[ht]
	\centering
	\includegraphics[height=0.20\textheight]{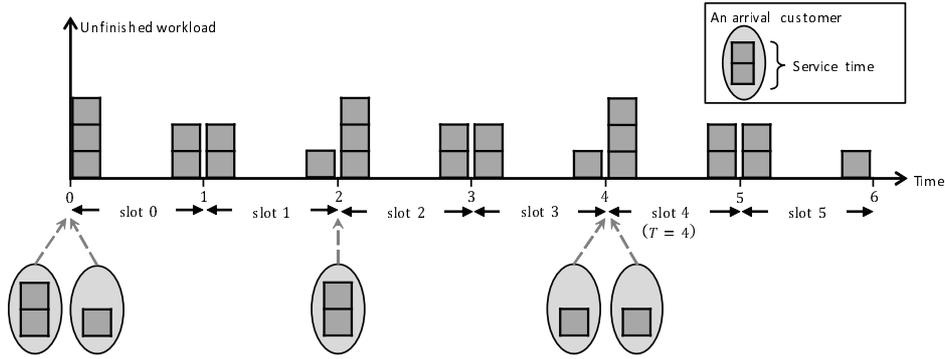}
	\caption{A sample path of the unfinished workload for $T = 4$.}
	\label{fig:samplepath}
\end{figure}

Throughout the paper, unless otherwise stated, Assumptions~\ref{assumpt-V_0-} and \ref{assumpt-service} are valid. Under these assumptions, a sample path of the (integer-valued) unfinished workload is given in Figure~\ref{fig:samplepath}.

We now present the definition of the queueing game studied in this paper. 
\begin{definition}[Poisson queueing game]
\label{def:Def_of_population}
{\it Potential} customers (i.e., candidates for players of the queueing game) are homogeneous. The number of them is countably infinite, and they are labeled with positive integer numbers $1,2,\dots$. With equally probability, each potential customer decides whether or not to arrive at the queueing system, independently of the other customers. Let $\bbA \subseteq \bbN$ denote the label set of the (actual) customers that decide to arrive at the queueing system, and assume that $|\bbA|$ follows a Poisson distribution with mean $\lambda$, i.e.,
\begin{eqnarray}
\PP( |\bbA| = n) = e^{-\lambda} {\lambda^n \over n!},
	\qquad n \in \bbZ_{+}.
\label{hm-add-170927-01}
\end{eqnarray}
The strategy of a customer is expressed by her/his arrival-time distribution, which is a probability distribution on $\bbT$. Let $\calP$ denote the set of all probability distributions on $\bbT$, and assume that customer $i$ chooses her/his arrival slot independently of the other customers, following a probability distribution $\vc{p}^{(i)}\in \calP$. The set $(\vc{p}^{(i)};i\in\bbA) \in \calP^{|\bbA|}$ is referred to as the {\it arrival-time distribution profile}, and the cost assessed by customer $i$, $i \in \bbA$, is denoted by $C^{(i)}$. Each customer chooses her/his arrival-time distribution with the goal of minimizing her/his cost. As a result, our queueing game is a {\it Poisson game} and is  represented by $(\bbA, \calP, \calC)$, where $\calC = (C^{(i)})_{i \in \bbA}$.
\end{definition}

\begin{remark}
The Poissonian population of customers is justified under a situation where the Poisson law of small numbers holds (see, e.g., \citealt[Exercise 4, page 181]{Chun01}). Furthermore, the Poissonian population has the property ``environmental equivalence" holds (see \citealt[Theorem~2]{Myer98}), that is,
\begin{eqnarray}
\PP( |\bbA^{(-i)}| = n \mid i \in \bbA) 
= \PP( |\bbA| = n)
= e^{-\lambda} {\lambda^n \over n!},
	\qquad n \in \bbZ_{+},
\label{hm-add-190517-01}
\end{eqnarray}
where $\bbA^{(-i)} = \bbA \setminus \{i\}$.
Thus, 
the property ``environmental equivalence" implies that the posterior belief of each (actual, not potential) customer about the number of the other ones is identical to the distribution of the total number of customers observed from the exterior.
\end{remark}

As described above, potential and thus actual customers (joining the system) are homogeneous. We thus focus on a symmetric (mixed-strategy) Nash equilibrium, where all the customers choose a common arrival-time distribution such that their expected waiting times are identical and such that they have no incentive to change its arrival-time distribution. To discuss such an equilibrium, we assume the following.
\begin{assumption}[Assessment for the behavior of the other customers]
\label{def:utility-payoff}
Each customer $i \in \bbA$ expects that the other customers in $\bbA^{(-i)}$ follow a single arrival-time distribution $\vc{p} \in \calP$, i.e., $\vc{p}^{(j)} = \vc{p} \in \calP$ for all $j \in \bbA^{(-i)}$. Consider a situation where customer $i$ chooses the arrival-time distribution $\vc{q}\in\calP$, i.e., $\vc{p}^{(i)} = \vc{q}:=(q_t ; t \in \bbT) \in \calP$ and assesses that the other customers in $\bbA^{(-i)}$ arrive with the arrival-time distribution $\vc{p}\in\calP$. Denote this situation by $(\vc{p}^{(i)} =\vc{q}, \vc{p})$. In situation $(\vc{p}^{(i)} =\vc{q}, \vc{p})$, the cost of customer $i$ is given by
\begin{equation*}
C^{(i)}(\vc{p}^{(i)} =\vc{q}, \vc{p}) := \sum_{t \in \bbT} q_{t} w_t(\vc{p}),
\qquad i \in \bbA,
\end{equation*}
where $w_t(\vc{p})$ denotes the expected waiting time of customer $i$  in situation $(\vc{p}^{(i)} =\vc{q}, \vc{p})$ provided that she/he is a member of the customers arriving in slot $t$. Since all customers are homogeneous, we simply write $C(\vc{q}, \vc{p})$ for $C^{(i)}(\vc{p}^{(i)} =\vc{q}, \vc{p})$ and thus 
\begin{equation}
C(\vc{q}, \vc{p}) = \sum_{t \in \bbT} q_{t} w_t(\vc{p}).
\label{defn-w(q|p)}
\end{equation}
\end{assumption}

We now provide the strict definition of a symmetric Nash equilibrium in our Poisson queueing game.
\begin{definition}[A symmetric Nash equilibrium]
An arrival-time distribution $\vc{p}^{*}:=(p_t^*;t\in\bbT)$ is said to be an {\it equilibrium arrival-time distribution} if and only if
\begin{align}
\label{eqn:Def_NE} 
C(\vc{p}^{*}, \vc{p}^{*}) \le C(\vc{q}, \vc{p}^{*}), \quad \forall \vc{q} \in \calP.
\end{align}
Furthermore, an arrival-time distribution profile $(\vc{p}^{(i)})_{i \in \bbA}$ is a symmetric Nash equilibrium if $\vc{p}^{(i)} = \vc{p}^{*}$ for all $i \in \bbA$.
\end{definition}

The following lemma presents an equivalent condition for (\ref{eqn:Def_NE}). 
\begin{lemma}
\label{defn:symmetric_Nash_equilibrium}
An arrival-time distribution $\vc{p}^{*} := (p^{*}_{t} ; t \in \bbT) \in \calP$ is the {\it equilibrium arrival-time distribution} if only if the following
holds: 
\begin{align} 
\label{eqn:Equilibrium_condition}
p^{*}_{t} > 0 \iff \ w_t(\vc{p}^{*}) = \min_{u \in \bbT} w_{u}(\vc{p}^{*}) =: w^*.
\end{align}
The value $w^{\ast}$ is called an {\it equilibrium expected waiting time}.
\end{lemma}
\proof We first prove the sufficiency of the statement. If (\ref{eqn:Equilibrium_condition}) holds, then, for any $\vc{q} \in \calP$, 
\begin{eqnarray*}
C(\vc{q}, \vc{p}^{*})
&=& \sum_{t \in \bbT} q_{t} w_t(\vc{p}^{*})
\nonumber
\\
&\ge&   \sum_{t \in \bbT} q_{t} \cdot w^*
= \sum_{t \in \bbT;\, p^{*}_{t}>0} p^{*}_{t} \cdot w^*
= \sum_{t \in \bbT} p^{*}_{t} w_t(\vc{p}^{*})
\nonumber
\\
&=& 
C(\vc{p}^{*}, \vc{p}^{*}),
\end{eqnarray*}
which shows that (\ref{eqn:Def_NE}) is satisfied.

Next we prove the necessity of the statement. 
Assume that (\ref{eqn:Def_NE}) holds. For proof by contradiction, suppose that the negation of (\ref{eqn:Equilibrium_condition}) holds, that is, there exist some $\sigma,u \in \bbT$ such that $p^{*}_{\sigma} > 0$, $p^{*}_{u} > 0$, and $w_{\sigma}(\vc{p}^{*}) > w_{u}(\vc{p}^{*})=w^*$. Fix $\vc{q}=(q_t;t\in\bbT)$ such that 
\begin{eqnarray*}
q_{\sigma} &=& 0, \quad q_{u} = p^{*}_{\sigma} + p^{*}_{u},
\\
q_{t} &=& p^{*}_{t},\quad t \neq \sigma,u. 
\end{eqnarray*}
We then have
\begin{eqnarray*}
C(\vc{q}, \vc{p}^{*})
&=& (p^{*}_{\sigma} + p^{*}_{u}) w_u(\vc{p}^{*}) 
+ \sum_{t\in\bbT\setminus \{\sigma,u\}} p^{*}_{t}  w_t(\vc{p}^{*})
\nonumber
\\
&<& p^{*}_{\sigma} w_{\sigma}(\vc{p}^{*}) + p^{*}_{u} w_u(\vc{p}^{*}) 
+ \sum_{t\in\bbT\setminus\{\sigma,u\}} p^{*}_{t} w_t(\vc{p}^{*})
\nonumber
\\
&=& \sum_{t\in\bbT} p^{*}_t w_t(\vc{p}^{*})
= C(\vc{p}^{*}, \vc{p}^{*}),
\end{eqnarray*}
which contradicts (\ref{eqn:Def_NE}). Therefore, (\ref{eqn:Equilibrium_condition}) is satisfied.
$\square$

\begin{remark}
The condition (\ref{eqn:Equilibrium_condition}) is equivalent to the one considered in Section~6 of \cite{Myer98}. 
\end{remark}

\section{The expected waiting time}
\label{sec:Expected waiting time}

In this section, we discuss the expected waiting time $w_t(\vc{p})$ of an arbitrary customer $i \in \bbA$ under Assumption~\ref{def:utility-payoff}. 
To this end, we consider a queueing process which results from the behavior of the customers in $\bbA^{(-i)}$, who follow a common arrival-time distribution $\vc{p}$.

Let $A_{t}$, $t \in \bbT$, denote the number of arrivals in slot $t$. 
According to Definition~\ref{def:Def_of_population} and Assumption~\ref{def:utility-payoff}, all the customers in $\bbA^{(-i)}$ choose their arrival slots independently one another, following a common arrival-time distribution $\vc{p}=(p_t;t\in\bbT) \in \calP$. This together with (\ref{hm-add-190517-01}) implies that
\begin{equation}
\PP(A_t = n)
= e^{-\lambda p_t} {(\lambda p_t)^n \over n!},
\qquad t \in \bbT,\ n \in \bbZ_+.
\label{eqn:Dist_of_A_{0}}
\end{equation}

Let $V_{t-}:=V_{t-}(\vc{p})$, $t \in \bbT$, denote the unfinished workload in the system immediately before slot $t$ starts (see Figure~\ref{fig:samplepath}),  i.e., at time $t-$. Let
\begin{equation}
v_t(k):=v_t(k; \vc{p})
= \PP(V_{t-} = k), \qquad t\in \bbT,~k\in \bbZ_+.
\label{defn-v_t(t|p)}
\end{equation}
In what follows, we sometimes use $v_t(k; \vc{p})$ instead of $v_t(k)$ to emphasize arrival-time distribution $\vc{p}$. Recall here that the system is empty immediately before slot 0 starts, which implies that $V_{0-} = 0$ and thus
\begin{equation}
\label{eqn:v_0(k)}
v_0(k)=
\left\{
\begin{array}{ll}
1, & k=0,
\\
0, & k\in\bbN.
\end{array}
\right.
\end{equation}

Let $S_t:=S_t(\vc{p})$, $t \in \bbT$, denote the total service time of the customers arriving in slot $t$, i.e.,
\begin{equation}
\label{eqn:Def_of_S_t}
S_t = \sum_{n=1}^{A_t} B_{t,n}, \qquad t\in \bbT,
\end{equation}
where the $B_{t,n}$, $n \in \{1,2,\dots,A_t\}$, denote
the service times of customers arriving in slot $t$. Recall that the number of arriving customers and their service times are independent. Thus, 
\begin{align} 
\EE[S_t] = \lambda p_t \beta,\qquad t\in \bbT.
\label{eqn-E[S_t]}
\end{align}
In addition, since the unfinished workload (if any) decreases by one in the end of each slot (see Assumption~\ref{assumpt-V_0-}), we have
\begin{align}
&&&&
V_{t-} &= \left( V_{(t-1)-} + S_{t-1} - 1 \right)^+,
& t & \in \bbT \setminus \{0\},&&&&
\label{recursion-V_t-}
\end{align}
where $(x)^+ = \max(x,0)$ for $x \in (-\infty,\infty)$. 

Let $(s_t(k); k \in \bbZ_+)$ denote the probability
mass function of $S_t$, i.e.,
\begin{align}
\label{eqn:s_t(k)}
s_t(k) = \PP(S_t=k), \quad t\in \bbT,~k \in \bbZ_+.
\end{align}
Note that $S_t$ follows a compound Poisson distribution from (\ref{eqn:Dist_of_A_{0}}), (\ref{eqn:Def_of_S_t}) and Assumption~\ref{assumpt-service}. Therefore, 
\begin{equation}
s_t(k)
= \sum_{n=0}^{\infty} e^{-\lambda p_t} {(\lambda p_t)^n \over n!} b^{*n}(k),
\qquad k\in\bbZ_+,
\label{eqn-s_t(k)}
\end{equation}
where $(b^{*n}(k);k\in\bbZ_+)$ denotes the $n$th-fold convolution of $(b(k);k\in\bbZ_+)$ itself. Note also that $(s_t(k);k\in\bbZ_+)$ in (\ref{eqn-s_t(k)}) is recursively computed as follows (see e.g., \citealt{Zhan16}): For each $t \in \bbT$,
\begin{align} 
s_{t}(0) &= e^{- \lambda p_{t}}, \label{eqn:s_{t}(0)}\\
s_{t}(k) &= \frac{\lambda p_{t}}{k} \sum_{m=1}^{k} m b(m) s_{t} (k-m), \quad k \in \bbN. 
\end{align}
Furthermore, (\ref{recursion-V_t-}) yields
\begin{eqnarray}
\label{eqn:v_t_and_s_t}
v_t(k)
&=& \left\{
\begin{array}{ll}
v_{t-1}(0) \{ s_{t-1}(0) + s_{t-1}(1) \}
+
v_{t-1}(1)s_{t-1}(0), & k=0,
\\
\dm\sum_{\ell=0}^{k+1}v_{t-1}(\ell)s_{t-1}(k+1-\ell), & k\in\bbN.
\end{array}
\right.
\end{eqnarray}

We now show the result on $(w_t(\vc{p}) ; t \in \bbT)$.
\begin{lemma}\label{lem-E[W_t]}
For any $\vc{p} \in \calP$,
\begin{align}
\label{eqn:ol{w}_{t}_thm-01}
 w_0(\vc{p}) &= \frac{\lambda p_{0}}{2} \beta,\\
\label{eqn:ol{w}_{t}_thm-02}
 w_t(\vc{p}) &= \sum_{j=0}^{t-1} \{ \lambda p_j \beta - 1 + e^{-\lambda p_j} v_j(0) \} + \frac{\lambda p_{t}}{2} \beta, & t \in \bbT \setminus \{0\}, 
\end{align}
%
where $\beta$ is the mean service time and $(v_t(k); t \in \bbT \setminus \{0\}, k\in\bbZ_+)$ is recursively determined by (\ref{eqn:v_t_and_s_t}).
\end{lemma}

\proof
Let $F_{t}:=F_t(\vc{p})$, $t\in\bbT$, denote the number of the customers
who arrive in slot $t$ and receive their services before customer $i$, provided that she/he is a member of the customers arriving in slot $t$. Since the customers arriving in one slot are served in random order, we have $\EE[F_{t}] = \lambda p_{t}/2$. It then follows that
\begin{align}
w_t(\vc{p}) 
&= \EE[V_{t-}] + \beta \EE[F_{t}]
\nonumber
\\
&= \EE[V_{t-}] + \frac{\lambda p_{t}}{2} \beta, \quad t \in \bbT.
\label{eqn-w_t(p)-02}
\end{align}
Since $V_{0-} = 0$, we have $\EE[V_{0-}] = 0$. Combining this and (\ref{eqn-w_t(p)-02}) yields (\ref{eqn:ol{w}_{t}_thm-01}).

In what follows, we prove (\ref{eqn:ol{w}_{t}_thm-02}). For this purpose, it suffices to show that 
\begin{align}
\EE[V_{t-}] 
= \sum_{j=0}^{t-1} 
\{ \lambda p_j \beta - 1 + e^{-\lambda p_j} v_j(0) \}, \quad t \in \bbT.
\label{eqn-E[V_{t-}]}
\end{align}
It follows from (\ref{eqn-E[S_t]}) and (\ref{recursion-V_t-}) that, for $t \in \bbT \setminus \{0\}$,
\begin{eqnarray}
\EE[V_{t-}]
&=& \EE[( V_{(t-1)-} + S_{t-1} - 1 ) \one(V_{(t-1)-} + S_{t-1} \ge 1)  ] \nonumber
\\
&=& \EE[( V_{(t-1)-} + S_{t-1}) \one(V_{(t-1)-} + S_{t-1} \ge 1)  ] 
- \PP(V_{(t-1)-} + S_{t-1} \ge 1)
\nonumber
\\
&=& \EE[V_{(t-1)-}] + \EE[S_{t-1}] 
- \{ 1- \PP(V_{(t-1)-} + S_{t-1} = 0) \}
\nonumber
\\
&=& \EE[V_{(t-1)-}] + \lambda p_{t-1} \beta
- 1 + \PP(V_{(t-1)-} + S_{t-1} = 0).
\label{eqn:V_{t-}_second}
\end{eqnarray}
Note here that $V_{(t-1)-}$ and $S_{t-1}$ are independent and nonnegative random
variables. Thus, using (\ref{eqn:s_{t}(0)}), we obtain
\begin{eqnarray}
\PP(V_{(t-1)-} + S_{t-1} = 0) 
&=& \PP(S_{t-1} = 0)\PP(V_{(t-1)-}=0) \nonumber \\
&=& e^{-\lambda p_{t-1}}v_{t-1}(0). \label{eqn:independent_V_{(t-1)-}_S_{t-1}}
\end{eqnarray}
Substituting (\ref{eqn:independent_V_{(t-1)-}_S_{t-1}}) into
(\ref{eqn:V_{t-}_second}) yields
\[
\EE[V_{t-}]
= \EE[V_{(t-1)-}] + \lambda p_{t-1} \beta
- 1 + e^{-\lambda p_{t-1}}v_{t-1}(0),
\quad t \in \bbT \setminus \{0\},
\]
which implies that (\ref{eqn-E[V_{t-}]}) holds. 
$\square$

\medskip

\section{Equilibrium arrival-time distributions}\label{sec-arrival-distribution}

This section discusses equilibrium arrival-time distributions. We first describe the construction of an equilibrium arrival-time distribution. We then provide some numerical examples.

\subsection{Construction of an equilibrium arrival-time distribution}
\label{sec:Arrival-Time Dist for Equilibrium}

The following theorem is a fundamental result on equilibrium arrival-time distributions.
\begin{theorem}
	\label{Main-Theorem}
Statements (a) and (b) below hold. \hfill
	\begin{enumerate}
		\item[(a)] $p_{0}^{\ast} > 0$.
		\item[(b)] An equilibrium arrival-time distribution $\vc{p}^*=(p_{t}^{\ast};t \in \bbT)$ is a solution of the equation (\ref{eqn-x=f(x)}), or equivalently, $\vc{p}^*$ satisfies the recursion:
		\begin{align}
		p_{t}^{*}
		= \left( p_{0}^{*}
		- 2\sum_{j=0}^{t-1} 
		\left\{ 
		p_j^{*} + {- 1 + e^{-\lambda p_j^{*}} v_j(0; \vc{p}^*) \over \lambda \beta } 
		\right\} \right)^+ , \quad t \in \bbT \setminus \{0\}.
		\label{p_t*-1}
		\end{align}
	\end{enumerate}
\end{theorem}

\proof
We first prove statement (a) by contradiction. To this end, we suppose that $p_0^{\ast}=0$. It then follows from (\ref{eqn:ol{w}_{t}_thm-01}) and (\ref{eqn-w_t(p)-02}) that
\begin{eqnarray*}
0 =w_0(\vc{p}^*) &\ge& \min_{t\in\bbT}w_t(\vc{p}^*) \ge 0,
\\
w_t(\vc{p}^*) &\ge& {\lambda p_t^{\ast} \over 2} \beta,
\qquad \forall t \in \bbT.
\end{eqnarray*}
These imply that if there exists some $j\in \bbT\setminus\{0\}$ such that $p_j^*>0$ then
\[
w_j(\vc{p}^*) \ge {\lambda p_j^{\ast} \over 2} \beta > 0
= \min_{t\in\bbT}w_t(\vc{p}^*),
\]
which contradicts Lemma~\ref{defn:symmetric_Nash_equilibrium}. Therefore, $p_0^* > 0$.

Next we prove statement (b). It follows from $p_0^*>0$ and Lemmas~\ref{defn:symmetric_Nash_equilibrium} and \ref{lem-E[W_t]} that
\begin{equation}
w^* = {\lambda p_{0}^{*} \over 2}\beta,
\label{eqn-w^*=E[W_0]}
\end{equation}
and, for $t \in \bbT \setminus\{0\}$,
\begin{eqnarray}
w^* 
&=& \sum_{j=0}^{t-1} 
\{ \lambda p_j^* \beta - 1 + e^{-\lambda p_j^*} v_j(0; \vc{p}^*) \} 
+ {\lambda p_{t}^* \over 2} \beta,\qquad \mbox{if $p_{t}^{*} > 0$},
\label{hm-add-180525-01}
\\
w^* 
&<& \sum_{j=0}^{t-1} 
\{ \lambda p_j^* \beta - 1 + e^{-\lambda p_j^*} v_j(0; \vc{p}^*) \} 
+ {\lambda p_{t}^* \over 2} \beta,\qquad \mbox{if $p_{t}^{*} = 0$}.
\label{hm-add-180525-02}
\end{eqnarray}
Substituting (\ref{eqn-w^*=E[W_0]}) into (\ref{hm-add-180525-01}) and (\ref{hm-add-180525-02}) yields, for $t \in \bbT \setminus\{0\}$,
\begin{align*}
p_t^* &= {2 \over {\lambda \beta}}
\left(
{\lambda p_0^* \over 2} \beta
- \sum_{j=0}^{t-1} 
\{ \lambda p_j^{*} \beta - 1 + e^{-\lambda p_j^{*}} v_j(0; \vc{p}^*) \}
\right), &  &\mbox{if $p_{t}^{*} > 0$},
\\
{\lambda p_0^* \over 2} \beta
&- \sum_{j=0}^{t-1} \{ \lambda p_j^{*} \beta - 1 + e^{-\lambda p_j^{*}} v_j(0; \vc{p}^*) \} < 0,& &\mbox{if $p_{t}^{*} = 0$}.
\end{align*}
Combining these leads to (\ref{p_t*-1}). 
$\square$

\medskip

Theorem~\ref{Main-Theorem} shows that an equilibrium arrival-time distribution $\vc{p}^*$ can be computed if $p_0^*>0$ is given. However, we cannot identify an appropriate value of $p_0^*$ in advance. Therefore, we have to find it by trial and error. Based on this fact, we establish a procedure for constructing an equilibrium arrival-time distribution $\vc{p}^*$.

Let $\vc{x} = (x_t;t\in\bbT) \in \bbR_+^{T+1}$. For $t \in \bbT$, let 
$(s_t(k; \vc{x}); k \in \bbZ_+)$ denote a probability distribution such that 
\begin{equation}
s_t(k; \vc{x})
= \sum_{n=0}^{\infty} e^{-\lambda x_t} {(\lambda x_t)^n \over n!} b^{*n}(k),
\qquad k\in\bbZ_+.
\label{eqn-s_t(k|x)}
\end{equation}
Furthermore, for $t \in \bbT$, let $(v_t(k; \vc{x});k\in\bbZ_+)$ denote a probability distribution determined by the recursion (\ref{eqn:v_t_and_s_t}) with $(s_t(k); k \in \bbZ_+)$ replaced by $(s_t(k; \vc{x}); k \in \bbZ_+)$: For $t =1,2,\dots,T$,
\begin{eqnarray}
v_t(k; \vc{x})
&=& \left\{
\begin{array}{ll}
v_{t-1}(0; \vc{x}) \dm\sum_{\ell=0}^1 s_{t-1}(\ell; \vc{x}) 
+
v_{t-1}(1; \vc{x})s_{t-1}(0; \vc{x}), & k=0,
\\
\rule{0mm}{5mm}
\dm\sum_{\ell=0}^{k+1}v_{t-1}(\ell; \vc{x})s_{t-1}(k+1-\ell; \vc{x}), & k\in\bbN,
\end{array}
\right.
\label{defn-v_t(k|x)}
\end{eqnarray}
where
\begin{equation}
v_0(k; \vc{x})
= \left\{
\begin{array}{ll}
1, & k=0,
\\
0, & k\in\bbN.
\end{array}
\right.
\label{defn-v_0(k|x)}
\end{equation}

Let $\vc{f}(\vc{x}) = (f_0(\vc{x}),f_1(\vc{x}),\dots,f_T(\vc{x}))$ for $\vc{x} = (x_t;t\in\bbT) \in \bbR_+^{T+1}$, where $f_t$ denotes a function $\bbR_+^{T+1} \to \bbR_+$ such that
\begin{equation}
f_t(\vc{x})
=
\left( x_0
- 2\sum_{j=0}^{t-1} 
\left\{ 
x_j + {- 1 + e^{-\lambda x_j} v_j(0; \vc{x}) \over \lambda \beta } 
\right\} \right)^+,\qquad t \in \bbT.
\label{defn-f_t(x)}
\end{equation}
Note that $f_0(\vc{x}) = x_0$ because the empty sum is defined as zero. It then follows from Theorem~\ref{Main-Theorem} that an equilibrium arrival-time distribution $\vc{p}^*$ is a solution of the equation:
\begin{equation}
\vc{x} = \vc{f}(\vc{x}),\qquad \vc{x} \in \bbR_+^{T+1}.
\label{eqn-x=f(x)}
\end{equation}

We now define $\vc{x}^*:=(x_t^*;t\in\bbT)$ as a vector such that $x_0^* \ge 0$ and the $x_t^*$, $t \in \bbT \setminus \{0\}$, are recursively determined by
\begin{eqnarray}
x_t^*
=
\left( x_0^*
- 2\sum_{j=0}^{t-1} 
\left\{ 
x_j^* + {- 1 + e^{-\lambda x_j^*} v_j(0; \vc{x}^*) \over \lambda \beta } 
\right\} \right)^+,\quad t \in \bbT \setminus \{0\}.
\label{defn-x_t^*}
\end{eqnarray}
Clearly, $\vc{x}^*$ is a solution of the equation (\ref{eqn-x=f(x)}).
\begin{theorem}\label{Main-Theorem-02}
	Statements (a)--(c) below hold. \hfill
	\begin{enumerate}
		\item[(a)] $x_{0}^{\ast} = 0$ if and only if $\vc{x}^* = \vc{0}$.
		\item[(b)] For each $t \in \bbT \setminus \{0\}$, $x_{t}^{*}$ is a continuous function of $x_{0}^{\ast}$.
		\item[(c)] There exists a solution of the equation (\ref{eqn-x=f(x)}) with $\sum_{t=0}^Tx_t = 1$.
	\end{enumerate}
\end{theorem}

\proof
We first prove statement (a).  Since the if-part is obvious, we prove the only-if part. It follows from  (\ref{defn-v_0(k|x)}) and (\ref{defn-x_t^*}) that
\begin{equation}
x_1^*
=
\left( x_0^*
- 2
\left\{ 
x_0^* + {- 1 + e^{-\lambda x_0^*} \over \lambda \beta } 
\right\} \right)^+.
\label{eqn-x_1^*}
\end{equation}
Therefore, if $x_0^* = 0$, then $x_1^* = 0$. We now suppose that
\begin{equation}
x_0^* = x_1^* = \cdots = x_{t-1}^* = 0\quad 
\mbox{for some $t\in \bbT \setminus \{0\}$}.
\label{cond-x_0^*=0}
\end{equation}
It then follows from (\ref{eqn-s_t(k|x)}) that $s_j(0; \vc{x}^*) = 1$ for all $j \in \{0,1,\dots,t-1\}$. Using this and (\ref{defn-v_t(k|x)}), we have
\begin{eqnarray}
v_j(0; \vc{x}^*) = 1\quad \mbox{for all $j \in \{0,1,\dots,t\}$}.
\label{eqn-v_t(0|x^*)=1}
\end{eqnarray}
Substituting (\ref{cond-x_0^*=0}) and (\ref{eqn-v_t(0|x^*)=1}) into (\ref{defn-x_t^*}) yields $x_t^* = 0$. By induction, statement (a) holds.

Next we prove statement (b) by induction. Equation (\ref{eqn-x_1^*}) shows that $x_1^{*}$ is continuous with respect to $x_0^{*}$. 
Suppose that there exists some $t \in \{2,3,\dots,T-1\}$ such that, for each $j=1,2,\dots,t-1$, $x_j^{*}$ is continuous with respect to $x_0^{*}$. Note that, for all $\varepsilon > 0$ and $\ell \in \bbT$,
\begin{eqnarray*}
	\left|
	\sum_{n=0}^{\infty} 
	e^{-\lambda (x_{\ell}^*+\varepsilon)} 
	{\{\lambda (x_{\ell}^*+\varepsilon)\}^n \over n!}b^{*n}(k)
	- 
	\sum_{n=0}^{\infty} e^{-\lambda x_{\ell}^*}{(\lambda x_{\ell}^*)^n \over n!}  b^{*n}(k)
	\right|
	\le 2.
\end{eqnarray*}
Thus, by the dominated convergence theorem, we obtain
\begin{eqnarray*}
	\lim_{\varepsilon\to0}
	\left|
	\sum_{n=0}^{\infty} e^{-\lambda (x_{\ell}^*+\varepsilon)} 
	{\{\lambda (x_{\ell}^*+\varepsilon)\}^n \over n!}b^{*n}(k)
	- 
	\sum_{n=0}^{\infty} e^{-\lambda x_{\ell}^*}{(\lambda x_{\ell}^*)^n \over n!}  b^{*n}(k)
	\right|
	= 0.
\end{eqnarray*}
This together with (\ref{eqn-s_t(k|x)}) implies that, for any $\ell \in \bbT$, $s_{\ell}(k; \vc{x}^*)$ is continuous with respect to $x_{\ell}^*$. Therefore, it follows from (\ref{defn-v_t(k|x)}) and the assumption of induction that, for each $j=1,2,\dots,t-1$, $v_j(k; \vc{x}^*)$ is continuous with respect to $x_0^*$. As a result,  $x_t^*$ given in (\ref{defn-x_t^*}) is continuous with respect to $x_0^{*}$. 

Finally, we prove statement (c). Let
\begin{equation}
G(x_0^*) = \sum_{t=0}^T x_t^*,\qquad x_0^* \ge 0,
\label{defn-G(x_0^*)}
\end{equation}
which is continuous with respect to $x_0^{*}$ due to statement (b). Furthermore, 
\[
G(0) = 0,
\qquad
G(1) \ge 1.
\]
Therefore, there exists some $\alpha > 0$ such that $G(\alpha) = 1$.
The proof has been completed.
$\square$

Based on Theorem~\ref{Main-Theorem-02}, we propose an algorithm for computing an equilibrium arrival-time distribution $\vc{p}^{\ast}= (p_{t}^{\ast};t\in \bbT )$.

\begin{algorithm}[H]

	\caption{Computing an equilibrium arrival-time distribution}\label{alg:Algorithn_for_Equilibrium}
	{\bf Input}: Sufficient small $\varepsilon,\delta \in (0,1)$
	\\ 
	{\bf Output}: 
	Equilibrium arrival-time distribution $\vc{p}^{\ast} = (p_{t}^{\ast};
	t\in \bbT)$
	
	\medskip
	
	Set $k=1$ and iterate the following:
	
	
	\begin{enumerate}
		\setlength{\parskip}{1mm} 
		\setlength{\itemsep}{1mm} 
		\item[1.] Set $x_0^{\ast} = k\varepsilon$.
		\item[2.] For $t = 1,2, \ldots, T$, compute $x_t^{\ast}$
		by (\ref{defn-x_t^*}).
		\item[3.] If $| 1 - \sum_{t \in \bbT} x_{t}^{\ast} |
		< \delta$, then set $\vc{p}^{\ast} = \vc{x}^{\ast}$, and otherwise increment $k$ by one and go to Step 1.
	\end{enumerate}
\end{algorithm}
\begin{remark}
\label{rem:Uniqueness}
Theorem~\ref{Main-Theorem} simply ensures that there exists at least one solution of the equation (\ref{eqn-x=f(x)}), and it does not ensure that the equation has an exactly one solution. Unfortunately, we have tried but failed to prove the uniqueness of the equilibrium arrival-time distribution. Although, in Section 6, we provide some numerical results that imply the uniqueness of it, they cannot completely deny there exist more than one equilibrium arrival-time distribution. All we can say is that Algorithm~\ref{alg:Algorithn_for_Equilibrium} yields an equilibrium arrival-time distribution achieving the smallest expected waiting time even if other equilibrium arrival-time distributions exist.
\end{remark}

\subsection{Numerical examples}
\label{sec:Numerical_Examples_No_EA}

In this subsection, we present some numerical examples to illustrate the effect of the service time distribution on the arrival strategy of customers. 

We set $\lambda = 5$ and $T = 20$. We then consider the three service-time distributions with mean $\beta \in (0,\infty)$:
\begin{itemize}
	\item Case 1 (Deterministic distribution):
	\[
	b(k)
	= 
	\left\{
	\begin{array}{l@{~~~}l}
	1, & k = \beta,
	\\
	0, & k \neq \beta,
	\end{array}
	\right.
	\]
	where the coefficient of variation (CV) is equal to zero.
	\item Case 2 (Geometric distribution):
	\begin{eqnarray*}
		b(k) &=& 
		\left\{
		\begin{array}{l@{~~~}l}
			0, & k=0,
			\\
			\dm{1 \over \beta} \left(1 - {1 \over \beta}\right)^{k-1},& k \ge 1,
		\end{array}
		\right.
	\end{eqnarray*}
	where the CV is equal to $\sqrt{1 - \beta^{-1}}$.
	\item Case 3 (Mixtures of geometric distributions):
	\begin{eqnarray*}
		b(k) &=& 
		\left\{
		\begin{array}{l@{~~~}l}
			0, & k=0,
			\\
			p \dm{1 \over \beta_{1}} \left(1 - {1 \over \beta_{1}}\right)^{k-1} + 
			(1-p) \dm{1 \over \beta_{2}} \left(1 - {1 \over \beta_{2}}\right)^{k-1},& k \ge 1,
		\end{array}
		\right.
	\end{eqnarray*}
	where $0 < p < 1$ and $\beta_{i} \ge 1$ for $i = 1,2$ such that $p \beta_{1} + (1-p) \beta_{2} = \beta$.
\end{itemize}

We note that the CV of Case 3 is equal to
\[
\beta^{-1}\sqrt{ 2( p \beta_{1}^{2} + (1-p) \beta_{2}^{2} ) - \beta (1+\beta) },
\]
which implies that the parameters of Case 3 are not uniquely determined even if the value of the CV is given. To solve this problem, we fix the parameters such that $p \beta_{1} = 1$ and thus
\begin{align*} 
p &= 1 - \xi (\beta - 1), \quad 
\beta_{1} = p^{-1}, \quad
\beta_{2} = \xi^{-1},
\end{align*}
where
\begin{align*} 
\xi &= 
\frac{4(\beta - 1)}
{3 \beta (\beta - 1) + y^2 \beta^2 +
	\sqrt{\{ 3 \beta (\beta - 1) + y^{2} \beta^2 \}^2
		- 8 (\beta - 1)^2 \{ y^2 \beta^2 + \beta (\beta + 1) \} }}.
\end{align*}
This setting ensures that the CV is equal to $y$.

\begin{figure}[ht]
	\centering
	\includegraphics[height=0.25\textheight]{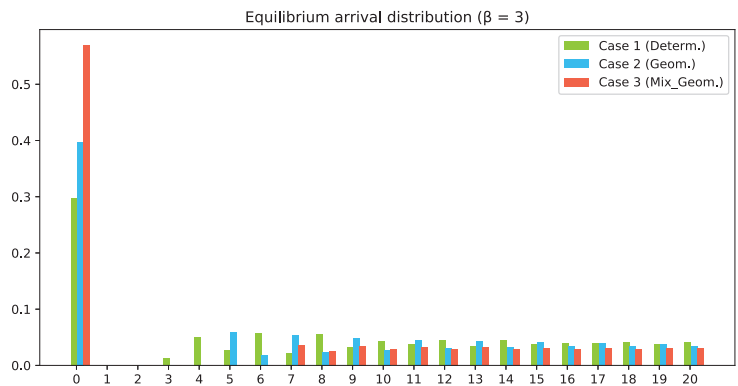}
	\caption{Equilibrium arrival-time distribution $(p_{0}^{*}, p_{1}^{*}, \ldots, p_{20}^{*})$ ($\beta = 3$)}
	\label{fig:graph_meanS_3}
\end{figure}

\begin{figure}[ht]
	\centering
	\includegraphics[height=0.25\textheight]{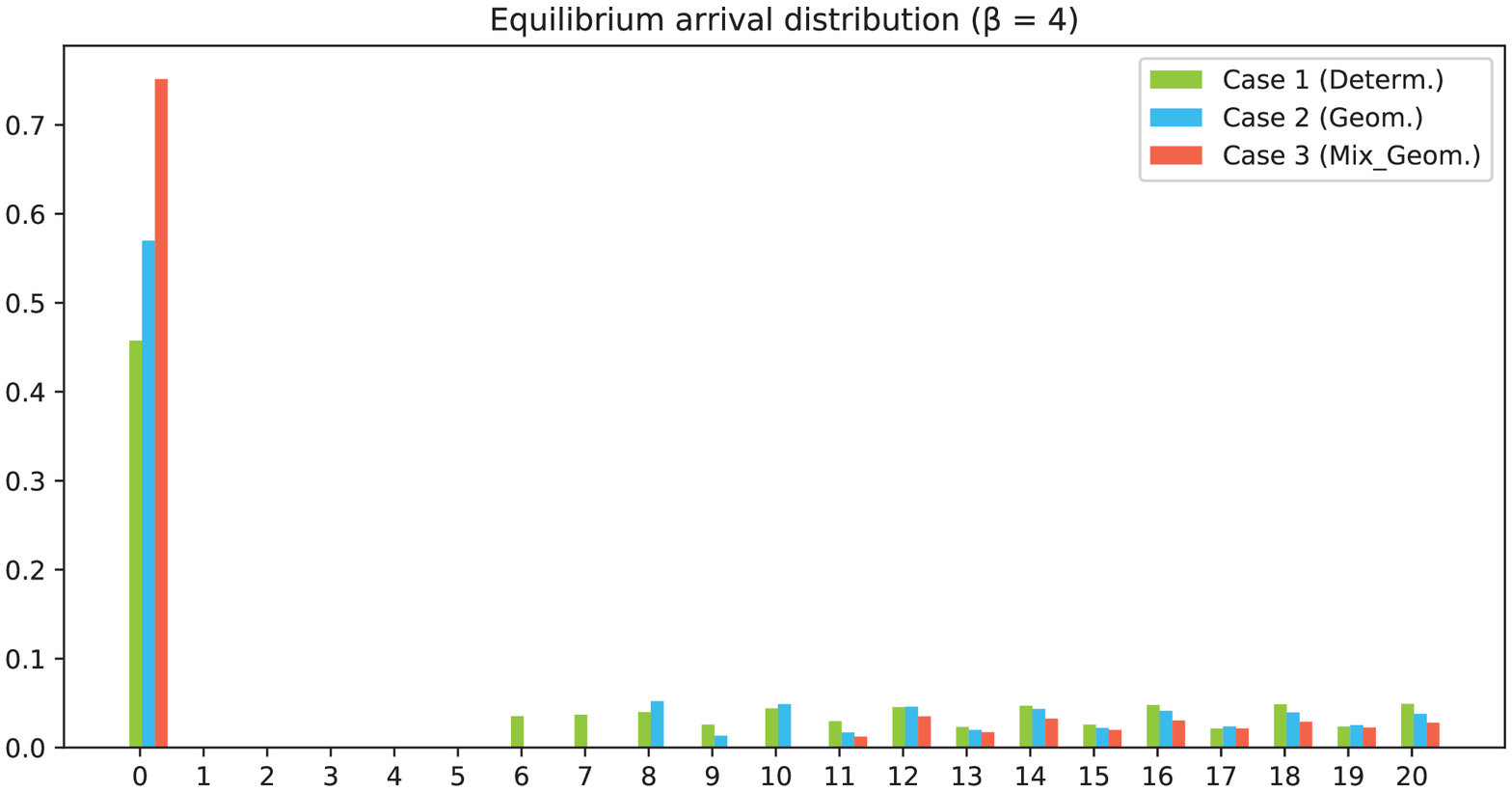}
	\caption{Equilibrium arrival-time distribution $(p_{0}^{*}, p_{1}^{*}, \ldots, p_{20}^{*})$ ($\beta = 4$)}
	\label{fig:graph_meanS_4}
\end{figure}

\begin{figure}[ht]
	\centering
	\includegraphics[height=0.25\textheight]{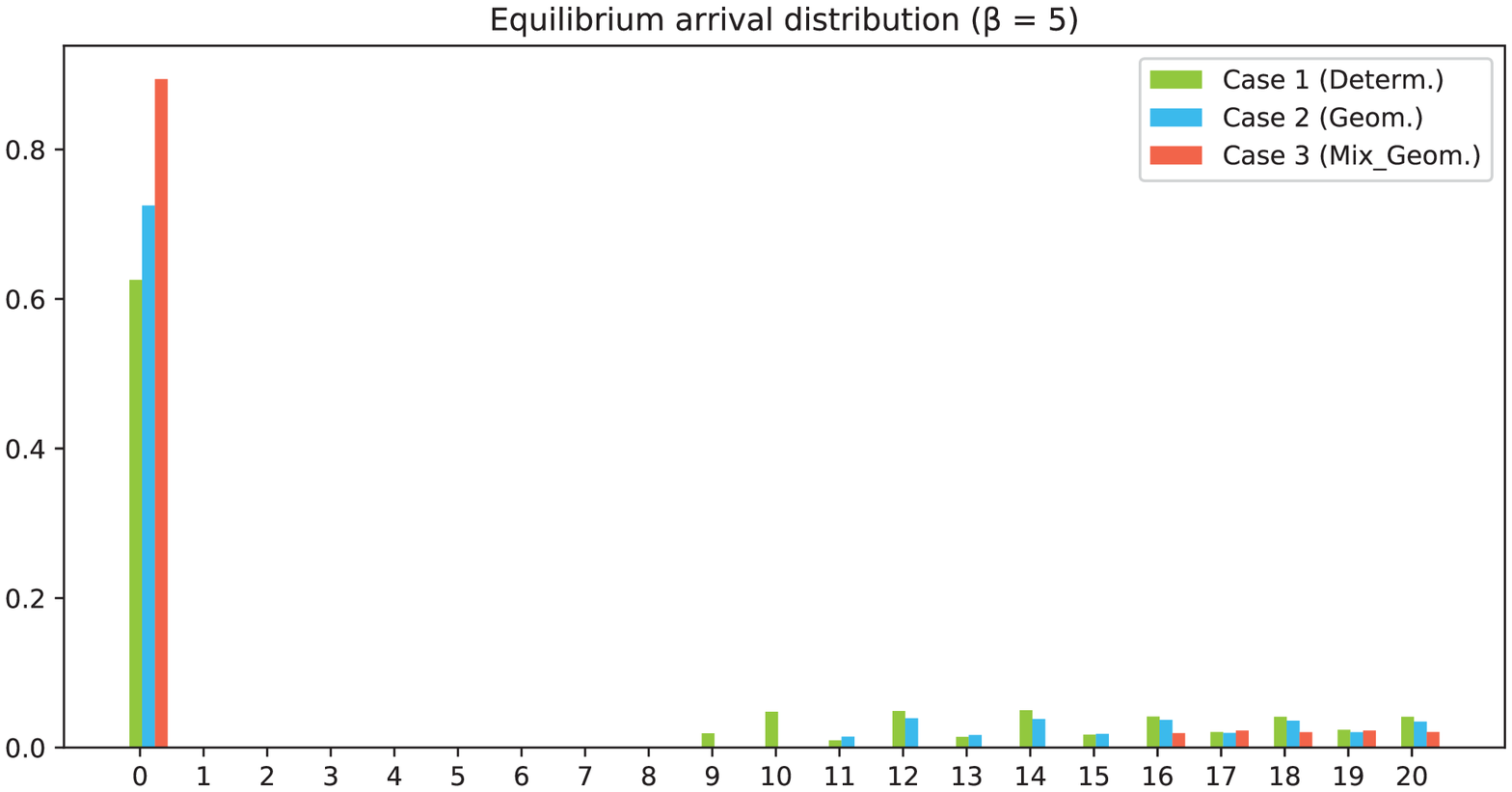}
	\caption{Equilibrium arrival-time distribution $(p_{0}^{*}, p_{1}^{*}, \ldots, p_{20}^{*})$ ($\beta = 5$)}
	\label{fig:graph_meanS_5}
\end{figure}

Using Algorithm~\ref{alg:Algorithn_for_Equilibrium}, we compute the equilibrium arrival-time distributions in Cases 1--3 for $\beta = 3, 4, 5$, which are shown in Figures~\ref{fig:graph_meanS_3}--\ref{fig:graph_meanS_5}. 
Furthermore, the equilibrium mean waiting times are summarized in Table~\ref{tab:table3_4_5}.
These results show that the probability $p_0^*$ of arriving in slot zero increases with the CV of the service-time distribution. 
This phenomenon is  explained as follows. Suppose that the CV of the service-time distribution is large and thus some customers demand large service times. Arriving after the arrival of such customers results in the large waiting time for service. Therefore, more customers rush to slot zero to minimize her/his waiting time.

%
\begin{table}[ht]
	\caption{Equilibrium expected waiting time for $\beta = 3, 4, 5$}
	\label{tab:table3_4_5}
	\begin{center}
		\begin{tabular}{|c|c|c|c|}
			\hline 
			& Case 1 & Case 2 & Case 3 \\ 
			\hline  \hline 
			$(\text{CV}, w^{*})$ for $\beta = 3$ & (0, 2.2)  & (0.82, 3.0) & (1.6, 4.3) \\ 
			\hline 
			$(\text{CV}, w^{*})$ for $\beta = 4$ & (0, 4.6)  & (0.87, 5.7) & (1.7, 7.5) \\ 			
			\hline 
			$(\text{CV}, w^{*})$ for $\beta = 5$ & (0, 7.8)  & (0.89, 9.1) & (1.8, 11.2) \\ 						
			\hline 
		\end{tabular}
	\end{center}
\end{table}

\begin{remark}
The shape of the equilibrium arrival-time distributions in Figures~\ref{fig:graph_meanS_3}-- \ref{fig:graph_meanS_5} are similar to ones in Figure 1 of \cite{Rapo04} as follows. In both of them, the opening slot has a mass of arrivals, then no arrivals occur in some subsequent slots, and arrivals resume and continue until the closing time.
\end{remark}

\section{Agent-based model}\label{sec:Simulation_Model}

We first propose an agent-based model, which simulates the behavior of customers that repeatedly join a single-server queue having an acceptance period and that try to minimize their waiting times in joining the queue according to their own histories of waiting times. We then compare the steady-state arrival-time distributions in our agent-based model with the equilibrium arrival-time distributions computed by Algorithm~\ref{alg:Algorithn_for_Equilibrium}.

\subsection{Model description}

Our agent-based model shares the basic assumptions with the queueing game considered in the previous sections. As with the queueing game, the agent-based model runs in discrete time, where a time is specified by ``day" and ``time slot (slot, for short)", for example, slot $t$ on the $\nu$-th day, where 
\[
t \in \bbT = \{0,1,\dots,T\},
\qquad \nu \in \bbN = \{1,2,\dots\}.
\]
Each day starts from slot zero and ends with the completion of serving the customers arriving on the day. The details of the model are described as follows.

The system has a single server and a waiting line of infinite capacity. On each day, the server accepts customers during slots $0,1,\dots,T$ (called the acceptance period), and the server continues to process the demands of customers in the system until it finishes processing them (even after slot $T$). Accepted customers are served on an FCFS basis, though the ones arriving in one slot are served in random order.

We assume that our agent-based model satisfies Assumptions~\ref{assumpt-V_0-} and \ref{assumpt-service}. Thus, the service times of customers are i.i.d.\ with the distribution $(b(k);k\in\bbN)$ with mean $\beta$. We also assume that there exist $N$ potential (homogeneous) customers of the system, which are labeled as the numbers in $\calN:=\{1,2,\dots,N\}$.
The number $N$ is fixed to be a sufficiently large constant.
On each day, the $N$ potential customers decide independently to arrive at the system or not. The probability of a customer arriving at the system is equal to $\Delta := \lambda / N$. 
Therefore, the number of arriving customers on each day follows a binomial distribution having mean $\lambda$. Since $N$ is sufficiently large, this binomial distribution is close to a Poisson distribution with mean $\lambda$, due to the Poisson law of small numbers (see e.g., \citealt[Exercise 4, page 181]{Chun01}).

To describe the strategy of customers, we introduce some definitions. 
Let $t_{\nu}^{(i)}$, $i \in \calN$, $\nu \in \bbN$, denote the arrival slot of customer $i$ on day $\nu$. 
Let $w_{\nu}^{(i)}$, $i \in \calN$, $\nu \in \bbN$, denote the waiting time of customer $i$ on day $\nu$.
If customer $i$ does not join the system on day $\nu$, set $t_{\nu}^{(i)} = \infty$ and $w_{\nu}^{(i)} = 0$. Moreover, for $i \in \calN$, $\nu \in \bbN$, and $t \in \bbT$, let
\begin{align*}
a_{\nu}^{(i)} = \sum_{t\in\bbT} a_{\nu,t}^{(i)}, \qquad
a_{\nu,t}^{(i)} = \sum_{k=1}^{\nu} \one(t_k^{(i)} = t),
\end{align*}
where $a_{\nu,t}^{(i)}$
denotes the total number of times that customer $i$ arrives in slot $t$ during the first $\nu$ days, and where $a_{\nu}^{(i)}$ denotes the total arrival days of customer $i$ during the first $\nu$ days. 
For convenience, let $a_{0}^{(i)} = 0$ for $i \in \calN$. 
Furthermore, for $i \in \calN$, $\nu \in \bbN$, and $t \in \bbT$, let
\begin{align*}
\ol{w}_{\nu,t}^{(i)} = {w_{\nu,t}^{(i)} \over \max(1, a_{\nu,t}^{(i)})}, \qquad
w_{\nu,t}^{(i)} = \sum_{k=1}^{\nu} \one(t_k^{(i)} = t) w_k^{(i)},
\end{align*}
where $w_{\nu,t}^{(i)}$ and $\ol{w}_{\nu,t}^{(i)}$ denote the sum and mean, respectively, of the waiting times that customer $i$ experiences when she/he arrives at slot $t$ during the first $\nu$ days.

As mentioned in the beginning of this section, every (arriving) customer tries to choose a slot that minimizes her/his waiting time for service according to her/his experience. Naturally, on the first day,  customers have no experiences, but they accumulate ``queueing experiences", i.e., their actual waiting times every time they join the system. Thus, every customer would gradually put weight on her/his own experience in choosing an arrival slot.
To express such a gradual shift toward experience-based decision, we introduce an increasing function $\theta: \bbR_+ \to [0,1]$ such that
\begin{equation}
\theta(0) = 0,\qquad \lim_{x \to \infty}\theta(x) = 1.
\label{cond-theta}
\end{equation}
We then assume that if customer $i$ ($i \in \calN$) arrives at the system on the $\nu$-th day ($\nu \in \bbN$) then the customer takes one of two options:
\begin{enumerate}
	\item[(i)] With probability $\theta(a_{\nu-1}^{(i)})$, to choose uniformly one slot $j \in \argmin_{t\in\bbT} \ol{w}_{\nu-1,t}^{(i)}$.
	\item[(ii)] With probability $1 - \theta(a_{\nu-1}^{(i)})$, to choose one slot uniformly from slots $0,1,\dots,T$.
\end{enumerate}
Note that, on the first day, customer $i$ ($i \in \calN$) chooses an arrival slot from slots $0,1,\dots,T$ with an equal probability, because of $a^{(i)}_{0} = 0$ and (\ref{cond-theta}).

Let $\vc{p}_{\nu}^{(i)}:=(p_{\nu,t}^{(i)};t\in\bbT)$, $i \in \calN$, $\nu \in \bbN$, denote the arrival-time distribution of customer $i$ on day $\nu$. Let
\[
\bbT_{\nu}^{(i)} = \argmin_{t \in \bbT} \ol{w}_{\nu,t}^{(i)},
\quad i \in \calN,\ \nu \in \bbN,
\]
which is the set of the slots that customer $i$ experiences the minimum mean waiting time during the first $\nu$ days.
It then follows from options (i) and (ii) that, for $i \in \calN$ and $t\in\bbT$,
\begin{align*}
p_{1,t}^{(i)} 
&= {1\over T + 1},
\\
p_{\nu,t}^{(i)}
&= \left\{
\begin{array}{ll}
\dm { 1 - \theta(a_{\nu-1}^{(i)}) \over T + 1}
+ 
{\theta(a_{\nu-1}^{(i)}) \over |\bbT_{\nu-1}^{(i)}|}, 
& t \in \bbT_{\nu-1}^{(i)},
\\
\dm { 1 - \theta(a_{\nu-1}^{(i)}) \over T + 1},
& t \not\in \bbT_{\nu-1}^{(i)},
\end{array}
\right.
\quad \nu \in \bbN \setminus \{1\}.
\end{align*}
Furthermore, let $\ol{\vc{p}}_{\nu}:=(\ol{p}_{\nu,t};t\in\bbT)$ and $\ol{w}_{\nu}$ , $\nu \in \bbN$, denote the averaged arrival-time distribution and waiting time, respectively, over all the customers on day $\nu$; that is,
\begin{align}
&&&&
\ol{p}_{\nu,t} &= {1 \over N} \sum_{i=1}^N p_{\nu,t}^{(i)},
& \nu &\in \bbN,~t\in\bbT,&&&&
\label{defn-ol{p}_{nu,t}}
\\
&&&&
\ol{w}_{\nu} &= {1 \over N} \sum_{i=1}^N {1 \over a_{\nu}^{(i)}} \sum_{k=0}^{\nu} w_k^{(i)},
& \nu &\in \bbN.&&&&
\label{defn-ol{w}_{nu}}
\end{align}

In the rest of this subsection, we explain our choice for function $\theta$.
Although there may be other reasonable options of function $\theta$, we fix
\begin{align}
\label{eqn:Def_SOHL} 
\theta(x) = 
\left\{
\begin{array}{ll}
0, & x = 0,
\\
\exp \left( \dm\frac{c_{1}}{1 - e^{ c_{2} x }} \right), & x > 0,
\end{array}
\right.
\end{align}
where $c_{1}$ and $c_{2}$ are positive parameters. The function $\theta$ in (\ref{eqn:Def_SOHL}) is a sigmoid-type function on the (nonnegative real) half line. Thus, for convenience, we refer to such a function as {\it sigmoid-on-the-half-line (SOHL for short) function}.

\begin{figure}[ht]
	\centering
	\includegraphics[width=\linewidth]{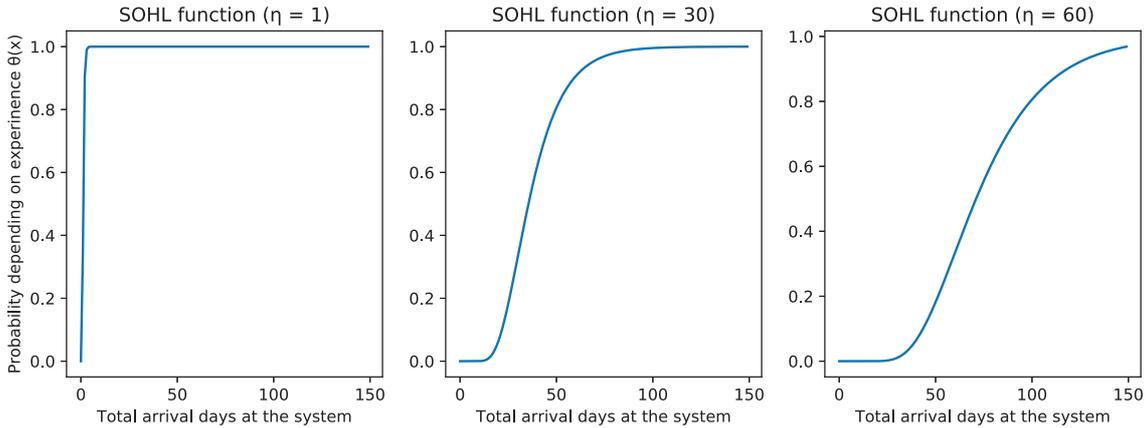}
	\caption{The shape of SOHL function with parameters in (\ref{eqn:Parameter_Settings_SOHL-c_1}) and (\ref{eqn:Parameter_Settings_SOHL-c_2})}
	\label{fig:shapesohl}
\end{figure}

By definition, the SOHL function $\theta$ has two parameters $c_1$ and $c_2$. To facilitate numerical experiments in the next subsection, we fix the parameters such that $\theta$ takes the value of $1/3$ at its unique {\it inflection point}. 
It follows from (\ref{eqn:Def_SOHL}) that $\theta$ is an  increasing twice-differentiable function satisfying (\ref{cond-theta}). Thus, there indeed exists a unique inflection point at $x=\eta$, given by
\begin{equation}
\eta = \frac{1}{c_{2}} \log \frac{c_{1} + \sqrt{c_{1}^{2} + 4}}{2}.
\end{equation}
Furthermore,
\begin{align} 
\theta(\eta) 
= \exp\left( \frac{2 c_{1}}{2 - c_{1} - \sqrt{c_{1}^{2} + 4}} \right).
\end{align}
Therefore, solving the equation $\theta(\eta) = 1/3$, we obtain 
\begin{equation}
c_{1} = \frac{(2 - \log 3) \log 3}{\log 3 - 1}.
\label{eqn:Parameter_Settings_SOHL-c_1} 
\end{equation}
Figure~\ref{fig:shapesohl} plots the SOHL functions $\theta$ with $\theta(\eta) = 1/3$ and $\eta = 1, 30,$ and $60$, where $c_1$ is given in (\ref{eqn:Parameter_Settings_SOHL-c_1}) and 
\begin{equation}
c_{2} = \frac{1}{\eta} \log \frac{c_{1} + \sqrt{c_{1}^2 + 4}}{2}.
\label{eqn:Parameter_Settings_SOHL-c_2} 
\end{equation}

\subsection{Numerical examples}

In this subsection, we present some numerical examples.
To this end, we set
\[
N = 100,\quad \calN = \{1,2,\ldots,N\},\quad T = 20,
\quad 
\lambda = 5.
\]
Thus, the number of potential customers is equal to $N=100$, and the acceptance period ends at slot $20$. The mean number of arriving customers on each day is equal to $\lambda = 5$, and therefore 
each customer joins the system with the probability $\Delta:=\lambda/N=5/100 = 0.05$. Theoretically, until day $\nu$, each customer has joined the system $\nu \times \Delta = 0.05 \nu$ times on average, in other words, 
\[
{1 \over 100}\sum_{i=1}^{100}a_{\nu}^{(i)} \approx  0.05 \nu.
\]
Furthermore, the parameters of the SOHL function are given by (\ref{eqn:Parameter_Settings_SOHL-c_1}) and (\ref{eqn:Parameter_Settings_SOHL-c_2}), which implies that $\theta(\eta) = 1/3$; that is, customers who have joined the system $\eta$ times chooses an arrival slot according to her/his experience with probability $1/3$.

\begin{figure}[ht]
	\centering
	\includegraphics[width=\linewidth]{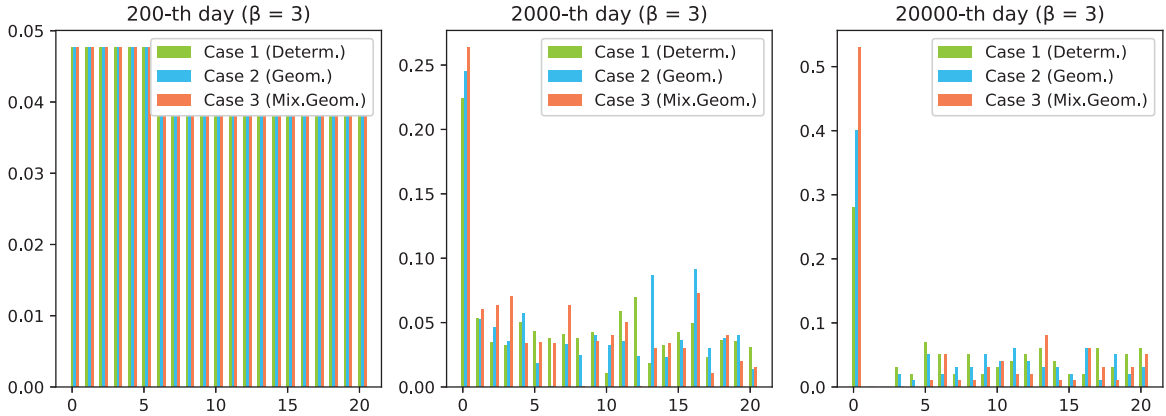}
	\caption{Averaged arrival-time distribution $\ol{\vc{p}}_{\nu}$ for $\nu=200, 2000, 2.0\times 10^4$ with $\beta = 3$ and $\eta = 30$.}
	\label{fig:Simu_graph_meanS_3}
\end{figure}

\begin{figure}[ht]
	\centering
	\includegraphics[width=\linewidth]{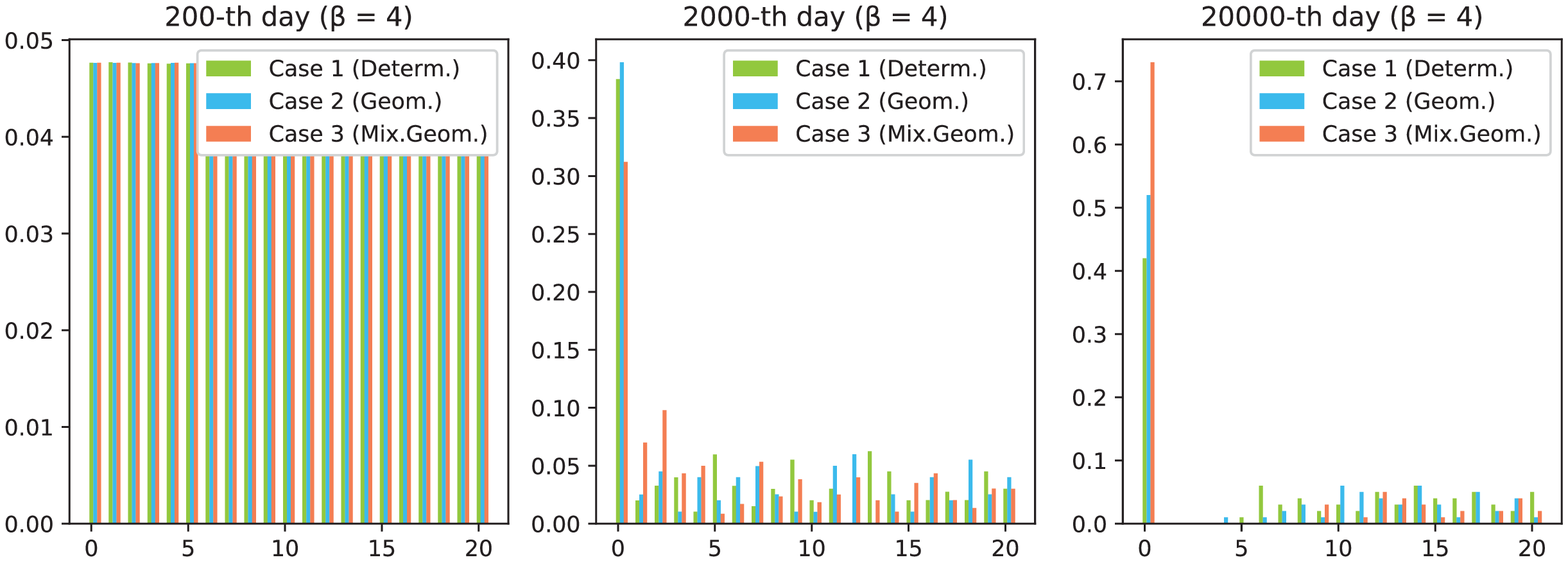}
	\caption{Averaged arrival-time distribution $\ol{\vc{p}}_{\nu}$ for $\nu=200, 2000, 2.0\times 10^4$ with $\beta = 4$ and $\eta = 30$.}
	\label{fig:Simu_graph_meanS_4}
\end{figure}

\begin{figure}[ht]
	\centering
	\includegraphics[width=\linewidth]{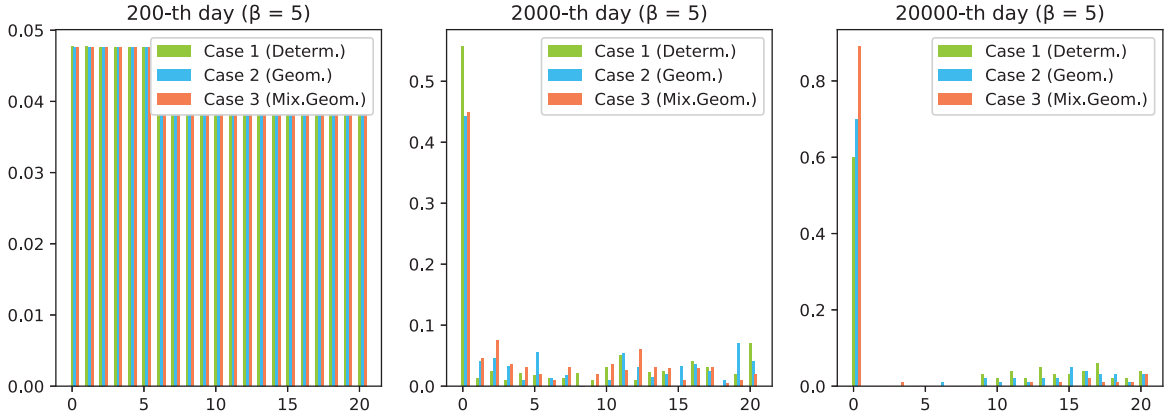}
	\caption{Averaged arrival-time distribution $\ol{\vc{p}}_{\nu}$ for $\nu=200, 2000, 2.0\times 10^4$ with $\beta = 5$ and $\eta = 30$.}
	\label{fig:Simu_graph_meanS_5}
\end{figure}

We first set $\eta = 30$.
Figure \ref{fig:Simu_graph_meanS_3} plots the averaged arrival-time distribution $\ol{\vc{p}}_{\nu}$ (defined in (\ref{defn-ol{p}_{nu,t}})) for $\nu= 200, 2000, 2.0\times 10^4$, where $\beta = 3$ ($\beta$ is the mean service time). Figures \ref{fig:Simu_graph_meanS_4} and \ref{fig:Simu_graph_meanS_5} also plot the results with $\beta = 4$ and $5$, respectively. 
These numerical results show that our agent-based model generates, in a sufficiently long run time, the averaged arrival-time distribution similar to the equilibrium arrival-time distribution (see Figures~\ref{fig:graph_meanS_3}--\ref{fig:graph_meanS_5}) computed by Algorithm~\ref{alg:Algorithn_for_Equilibrium}. The simulation results also adhere to the consequences of our mathematical model in that the probability of arriving in slot zero increases with $\beta$ and/or CV. As a result, we can say that our queueing game serves as a mathematical model of the real behavior of non-cooperative customers.

\begin{figure}[ht]
	\centering
	\includegraphics[width=\linewidth]{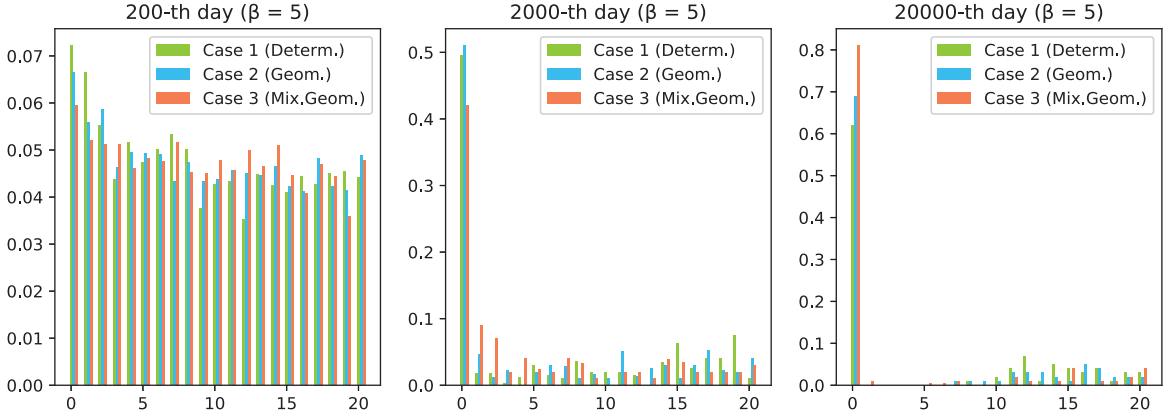}
	\caption{Averaged arrival-time distribution $\ol{\vc{p}}_{\nu}$ for $\nu=200, 2000, 2.0\times 10^4$ with $\beta = 5$ and $\eta = 1$.}
	\label{fig:Simu_graph_meanS_5_eta1}
\end{figure}

\begin{figure}[ht]
	\centering
	\includegraphics[width=\linewidth]{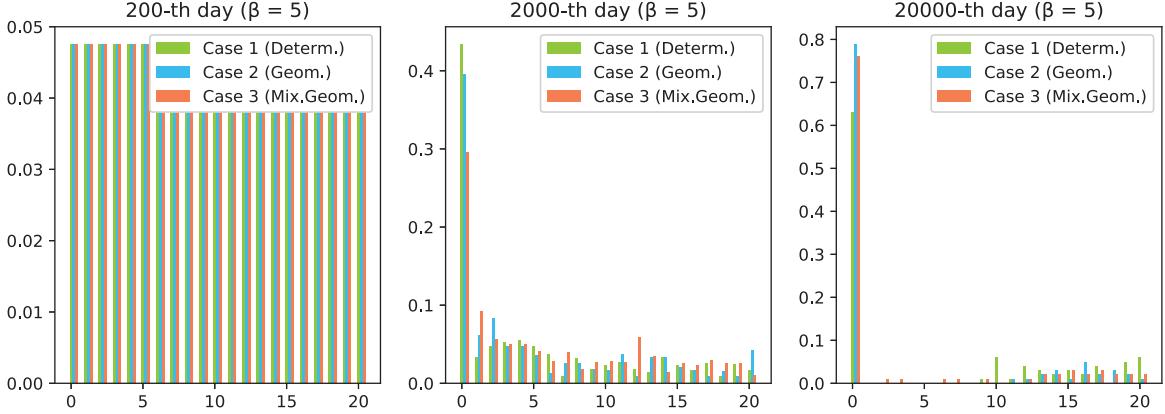}
	\caption{Averaged arrival-time distribution $\ol{\vc{p}}_{\nu}$ for $\nu=200, 2000, 2.0\times 10^4$ with $\beta = 5$ and $\eta = 60$.}
	\label{fig:Simu_graph_meanS_5_eta60}
\end{figure}

Next, we examine how the parameter $\eta$ of the SOHL function impacts on
	the averaged arrival-time distribution.
Figures \ref{fig:Simu_graph_meanS_5_eta1} and \ref{fig:Simu_graph_meanS_5_eta60} plot the averaged arrival-time distributions for $\nu= 200, 2000, 2.0\times 10^4$ and $\beta = 5$, where $\eta = 1$ and $\eta = 60$, respectively. In the case of $\eta = 1$ (Figure~\ref{fig:Simu_graph_meanS_5_eta1}), the customers choose their arrival slots according to their own experiences from very early stages, and thus the probability of arriving in slot zero gradually increases after only 200 days. In contrast, in the case of $\eta = 60$ (Figure~\ref{fig:Simu_graph_meanS_5_eta60}), the customers choose their arrival slots according to the uniform distribution with probability more than 20\% even 200 days have passed. For this reason, the probabilities of arriving in slots 1--3 in the case of $\eta=60$ (Figure~\ref{fig:Simu_graph_meanS_5_eta60}) are higher than those in the cases of $\eta = 1$ (Figure~\ref{fig:Simu_graph_meanS_5_eta1}) and $\eta=30$ (Figure~\ref{fig:Simu_graph_meanS_5}) even after the 2000 days have passed.
However, after a sufficiently long time, the averaged arrival-time distributions become quite similar to the equilibrium arrival-time distributions in Figure \ref{fig:graph_meanS_5} regardless of $\eta$.

\section{Concluding remarks}
\label{sec:Conclusion}

In this paper, we have studied a discrete-time single-server Poisson queueing game with a single acceptance period. For this queueing game, we have proved that there exists an equilibrium arrival-time distribution such that the expected waiting times of customers are identical.

\begin{figure}[ht]
	\centering
	\includegraphics[width=\linewidth]{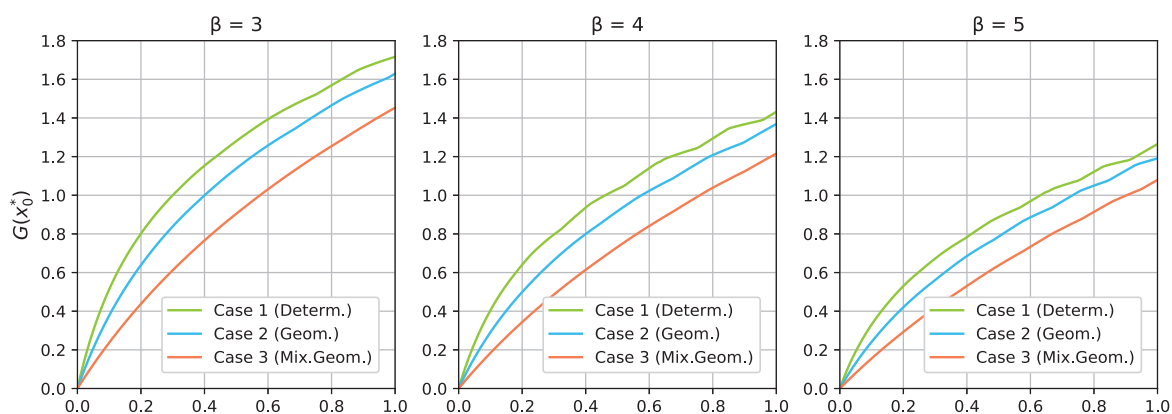}
	\caption{Monotonicity of $G(x_0^*)$ for three service-time distributions with $\beta = 3, 4, 5$.}
	\label{fig:sumgraph}
\end{figure}

We also have proposed an algorithm for computing an equilibrium arrival-time distribution achieving the smallest expected waiting time. Unfortunately, we have not proved that the equilibrium arrival-time distribution is unique, or equivalently, that there exists a unique solution of the equation (\ref{eqn-x=f(x)}) with $\sum_{t\in\bbT}x_t=1$. However, Figure~\ref{fig:sumgraph} shows the monotone increase of $G(x_0^*)$ with $x_0^*$ under the settings of Table~\ref{tab:table3_4_5}. These numerical results imply the uniqueness of the equilibrium arrival-time distribution.

%

We note that if the equilibrium arrival-time distribution is unique then it can be efficiently computed by the bisection method instead of Algorithm~\ref{alg:Algorithn_for_Equilibrium}. Therefore, it is a challenging and interesting problem to prove whether or not the equilibrium arrival-time distribution is unique.

%
%
%
%

Finally, we mention other two interesting directions for further research associated with this paper. One is to generalize our queueing model by taking into account earliness and tardiness costs to investigate the impact of these costs on the equilibrium arrival-time distribution. Assuming exponential service times, \cite{Sher17} considers a continuous-time Poisson queueing game with the linear combination of waiting, earliness, and tardiness costs. Discretizing the time axis would facilitate to handle the general service-time distribution in Poisson queueing games having such more flexible costs.

The other one is to study how the service discipline impacts on the equilibrium arrival-time distribution. \cite{Brei16} discuss the impact of FCFS, LCFS, and SIRO (Service In Random Order) on equilibria in a ``non-stochastic" queueing game with only three customers.  As far as we know, this service discipline problem has not yet considered in ``stochastic" queueing games having acceptance periods as with our queueing game, though it has done in ``To Join or Not to Join" models without acceptance periods (see, e.g., \citealt{Hass85}).


\section*{Acknowledgments}
The authors thank Mr.\ Hirotaka Kuwano for pointing out some minor errors in an early version of this paper. 
The authors also thank anonymous referees for their thoughtful comments.
The first, second, and third authors' works are supported by JSPS KAKENHI Grants  No.~JP16K21704, No.~JP18K11181, and No.~JP17K12980, respectively.

%
%
\bibliographystyle{elsarticle-harv}
%
%


\end{document}